%% file: main.tex
\documentclass[10pt]{amsart}
\input{macros.tex}
\begin{document}

\title[Smooth Siegel disks via semicontinuity]
{Smooth Siegel disks via semicontinuity: a remark on a proof of Buff
and Cheritat}

\author{Artur Avila}

\address{
Coll\`ege de France -- 3 Rue d'Ulm \\
75005 Paris -- France.
}
\email{avila@impa.br}

\thanks{Partially supported by Faperj and CNPq, Brazil.}   

\begin{abstract}

Recently, Xavier Buff and Arnaud Cheritat have provided an elegant proof
of the existence of quadratic Siegel disks with smooth boundary.
In this short note, we show how results of Yoccoz and Risler can be used to
conclude the same result.  Our proof is a small modification of the
argument given by Buff and Cheritat.

\end{abstract}

\setcounter{tocdepth}{1}

\maketitle


\input{smsie.tex}

\input{bib.tex}

\end{document}

%% file: macros.tex
\let\newpf\proof \let\proof\relax 
\newenvironment{pf}{\newpf[\proofname]}{\qed\endtrivlist}

\def\be{\begin{equation}}
\def\ee{\end{equation}}

\def\ba{{\begin{align}}}
\def\ea{{\end{align}}}

\def\0{{\mathbf 0}}

\newtheorem{thm}{Theorem}[section]
\newtheorem{cor}[thm]{Corollary}

\newtheorem{lemma}[thm]{Lemma}

\theoremstyle{remark}
\newtheorem{rem}{Remark}[section]

\numberwithin{equation}{section}

\def \bn {\hfill \\ \smallskip\noindent}

\theoremstyle{definition}

\def\proof{\bn {\bf Proof.} }

\def\note#1
{\marginpar
{\tiny $\leftarrow$
\par
\hfuzz=20pt \hbadness=9000 \hyphenpenalty=-100 \exhyphenpenalty=-100
\pretolerance=-1 \tolerance=9999 \doublehyphendemerits=-100000
\finalhyphendemerits=-100000 \baselineskip=6pt
#1}\hfuzz=1pt}

\newcommand{\dist}{\operatorname{dist}}

\newcommand{\C}{{\mathbb C}}
\newcommand{\D}{{\mathbb D}}

\newcommand{\Q}{{\mathbb Q}}
\newcommand{\R}{{\mathbb R}}

\def\B0{{\bold{0}}}

\catcode`\@=12

\def\Empty{}
\newcommand\oplabel[1]{
  \def\OpArg{#1} \ifx \OpArg\Empty {} \else
  	\label{#1}
  \fi}
		
\newcommand{\comm}[1]{}
\newcommand{\comment}[1]{}

%% file: smsie.tex
\section{Introduction}

Recently, in \cite {BC}, Xavier Buff and Arnaud Cheritat gave a new
proof of the following unpublished result of Perez-Marco: {\it
there exists a quadratic map with a Siegel disk whose boundary is a smooth
($C^\infty$) Jordan curve}.  Their proof involves both techniques of
renormalization of \cite {yoccozast} and estimates for parabolic explosion.

Our aim in this note is to show that the same result follows easily from
renormalization theory via two known results (of Yoccoz and Risler) by some
general abstract reasoning (which is really just a
small modification of \cite {BC}).

We would like to note that the method of parabolic explosion (coupled with
renormalization) allows much
greater control of the dynamics.  In particular, in \cite {BC} it is also
possible to conclude that the Siegel disks are accumulated by small cycles,
see also \cite {BC1} for even more dramatic applications.  However, we find
it worthwhile to investigate what is really needed for the argument, and we
hope that the present treatment could be used in situations which are
more general than the quadratic setting: indeed the proof works for
families of rational or entire maps which do not have non-Brjuno Siegel
disks, such as the families $z \to e^{2\pi i \alpha} z (1+z/d)^d$,
$d \geq 2$, and $z \to e^{2 \pi i \alpha} z e^z$ considered by Lukas Geyer
in \cite {G}.  We remark that it is conjectured that Siegel disks of
rational maps are always Brjuno.

In the last section we discuss the application of the method to the case of
Herman rings.  This application has been pointed out to us by Xavier Buff,
who had obtained this result earlier by other methods.

\section{Main result}

\subsection{Siegel disks}

Let $P_\alpha(z)=e^{2 \pi i \alpha} z+z^2$.
Let $r_\alpha$ be the conformal radius of the Siegel disk $\Delta_\alpha$ of
$P_\alpha$ if it exists, and let $r_\alpha=0$ otherwise.
If $r_\alpha>0$, let $L_\alpha:\D_{r_\alpha} \to \Delta_\alpha$ (where
$\D_r=\{|z|<r\}$)
be the uniformization map
satisfying $L_\alpha(0)=0$ and $DL_\alpha(0)=1$.
The function $L_\alpha$ satisfies the functional equation
$P_\alpha(L_\alpha(z))=L_\alpha(e^{2 \pi i \alpha} z)$.

Let $F_r$ be the space of holomorphic functions $f:\D_r \to \C$ with the
topology of uniform convergence on compact subsets of $\D_r$.
Let $E_r$ be a complete metric space of functions
$f:\D_r \to \C$.  We assume that for $r'>r$ we have $F_{r'} \subset E_r$
and that the inclusion is continuous.  For instance, $E_r$ can be taken as
the Fr\'echet space of $C^\infty$ functions $f:\overline \D_r \to \C$.
The requirements also allows one to consider (subspaces of)
certain spaces of quasianalytic functions, as the Banach space of
$C^\infty$ functions $f:\overline \D_r \to \C$ such that $f|\D_r$ is
holomorphic and $\sup_{r \geq 2}
\sup_{x \in \partial \D_r} \frac {|\partial^r f(x)|} {(r \ln r)^r}<\infty$.


\begin{thm} \label {2.1}

Let $r_{\alpha_0}>0$.  For every $\delta>0$, $0<r<r_{\alpha_0}$,
there exists $\alpha \in \R$ such that $|\alpha-\alpha_0|<\delta$,
$r_\alpha=r$, $L_\alpha|\D_r \in E_r$, and
$\dist_{E_r}(L_{\alpha_0}|\D_r,L_\alpha|\D_r)<\delta$.

\end{thm}

In order to prove Theorem \ref {2.1}, we will use properties of the function
$\alpha \to r_\alpha$.  Two of them are elementary:

(P1)\, $r_\alpha=0$ for a dense set of $\alpha$.  Indeed if
$\alpha=\frac {p} {q} \in \Q$ and $r_\alpha>0$ then $P_\alpha^q$ would have
to be the identity on a neighborhood of $0$, but $P_\alpha^q$ is a monic
polynomial of degree $2^q$.

(P2)\,
The function $\alpha \to r_\alpha$ is upper semicontinuous.  Indeed,
if $\alpha_n \to \alpha$ and $\inf_n r_{\alpha_n} \geq r'>0$
then Hurwitz Theorem $L_{\alpha_n}|\D_{r'}$ converges, in the topology
of $F_{r'}$ (and hence also in the topology of $E_r$, $r<r'$)
to a univalent fuction, which must coincide with $L_\alpha|\D_{r'}$.

We also need two non-elementary properties, which depend on renormalization
theory through results of Yoccoz and Risler.
Let us say that a function $h:\R \to \R$ is weakly lower
semicontinuous at $c \in \R$, we have
\be \label {alsc}
\min \{\limsup_{y \to c+} h(y),\limsup_{y \to c-} h(y)\} \geq h(c).
\ee

(P3)\, $r_\alpha$ is weakly lower semicontinuous when $\alpha$ is
non-Brjuno\footnote {Recall that $\alpha \in \R$ is called a Brjuno number
if it is irrational and $\sum \frac {\ln q_{n+1}} {q_n}<\infty$, where $q_n$
is the (increasing) sequence of denominators of the best rational
approximations of $\alpha$.}.  Indeed Yoccoz's Theorem \cite {yoccozast}
implies that
$r_\alpha=0$ for non-Brjuno numbers, so by (P1) and $r_\alpha \geq 0$,
$\alpha \in \R$, we see that $r_\alpha$ is even
continuous at non-Brjuno numbers.\footnote {This is the only place we
shall use special properties of the quadratic family.}

(P4)\, $r_\alpha$ is weakly lower semicontinuous when $\alpha$ is
Brjuno.  Indeed, by a result of Risler, if $\alpha$ is Brjuno then there
exists a set $\mathcal {B}_s$ restricted to which the function
$\alpha \to r_\alpha$ is continuous (see Proposition 10 of \cite {risler}),
and $\alpha$ is a Lebesgue density point of $\mathcal {B}_s$
(see Proposition 1 of \cite {risler} for other properties of the sets
$\mathcal {B}_s$).\footnote{It would be interesting to investigate if the
estimates of Risler are enough to conclude that $r_\alpha$ is weakly lower
semicontinuous also at non-Brjuno $\alpha$, as this would remove the
necessity of the step (P3) and make the whole argument much more general.}

\comm{
Let us say that a function $h:\R \to \R$ is weakly lower
semicontinuous if for every $c \in \R$, we have
\be \label {alsc}
\min \{\limsup_{y \to c+} h(y),\limsup_{y \to c-} h(y)\} \geq h(c).
\ee

The following result is the key step which uses the renormalization theory
through results of Yoccoz and Risler.

\begin{lemma}

The function $\alpha \to r_\alpha$ is weakly lower semicontinuous.

\end{lemma}

\begin{pf}

Since $r_\alpha \geq 0$, $\alpha \in \R$, it is enough to check (\ref
{alsc}) when $r_\alpha>0$.  By Yoccoz's Theorem \cite {yoccozast},
this implies that $\alpha$
satisfies the Brjuno condition.  By a result of Risler (see Proposition 10
of \cite {risler}), (\ref {alsc}) is
satisfied for $\alpha$ Brjuno: indeed for any sequence $\alpha_n \to \alpha$
such that $\Phi(\alpha_n) \to \Phi(\alpha)$, we have $\liminf r_{\alpha_n}
\geq r_\alpha$, where $\alpha \to \Phi(\alpha)$ is the arithmetic Yoccoz
function.
\end{pf}
}

The properties (P2-4) will be exploited through the following:

\begin{lemma}

If $h:\R \to \R$ is upper semicontinuous and weakly lower semicontinuous at
every $c \in \R$,
then $h$ satisfies the Intermediate Value Theorem.

\end{lemma}

\begin{pf}

Let $a<b$ such that $h(a) \neq h(b)$.  To fix ideas, assume $h(a)<h(b)$. 
Let $h(a)<x<h(b)$ and let $c=\inf \{a \leq y \leq b,\, h(y) \geq x\}$.  By
upper semicontinuity, $h(c) \geq x$, so $a<c \leq b$.  If $h(c)>x$,
by (\ref {alsc}) there exists $a<y<c$ such that $h(y)>x$, contradicting
the definition of $c$.  Thus $h(c)=x$ as required.
\end{pf}

Together with (P1), this yields:

\begin{cor} \label {2.4}

If $r_\alpha>0$ then for every $0<r<r_\alpha$ and $\epsilon>0$ there exists
$\alpha' \in \R$ such that $|\alpha'-\alpha|<\epsilon$ and $r_{\alpha'}=r$.

\end{cor}

\noindent {\it Proof of Theorem \ref {2.1}.}

Let $\beta_i$, $\epsilon_i$ be defined inductively as follows.  Let
$\beta_0=\alpha_0$ and $\epsilon_0=\delta$.  Assuming $\beta_i$,
$\epsilon_i$ defined, let $\epsilon_{i+1}<\epsilon_i/10$ be such that
$r_\beta<r_{\beta_i}+2^{-i}$ whenever $|\beta-\beta_i|<\epsilon_{i+1}$ (this
is possible by upper semicontinuity).  Then let
$\beta_{i+1}$ be such that $|\beta_{i+1}-\beta_i|<\epsilon_{i+1}/10$,
$r_{\beta_{i+1}}=\frac {r+r_{\beta_i}} {2}$, and
$\dist_{E_r}(L_{\beta_{i+1}}|\D_r,L_{\beta_i}|\D_r)<\epsilon_{i+1}/10$
(this is possible by Corollary \ref {2.4} and the proof of (P2) above).
It is easy to check that $\alpha=\lim \beta_i$ has the required properties.

\begin{rem}

One easily checks that this proof yields a Cantor set of $\alpha$ satisfying
the conclusions of Theorem \ref {2.1}.

\end{rem}

\subsection{Herman rings}

Fix $a>3$ and let $Q_\lambda(z)=e^{2 \pi i \lambda} z^2 \frac {z+a}
{1+az}$.  Then
$Q_\lambda(z)$, $\lambda \in \R$, is a diffeomorphism of $S^1$, and we let
$\alpha_\lambda$ be its rotation number.  Let $r_\lambda$ be half the
moduli of the Herman ring $\Xi_\lambda$ of $Q_\lambda$ if it exists, and let
$r_\lambda=0$ otherwise.
If $r_\lambda>0$, let $T_\lambda:A_{r_\lambda} \to \Xi_\lambda$ (where
$A_r=\{-r<\ln |z|<r\}$) be
the uniformization map satisfying $T_\lambda(1)=1$ and $DT_\lambda(1)>0$.
The map $T_\lambda$ satisfies the functional equation $T_\lambda(e^{2 \pi
i \alpha_\lambda} z)=Q_\lambda(T_\lambda(z))$.

Let $F_r$ be the space of holomorphic functions $f:A_r \to \C$ with the
topology of uniform convergence on compact subsets of $A_r$.
Let $E_r$ be a complete metric space of functions $f:A_r \to \C$.
For $r'>r$, we assume that $F_{r'} \subset E_r$ and that the inclusion
is continuous.

\begin{thm}

Let $r_{\lambda_0}>0$.  For every $\delta>0$, $0<r<r_{\lambda_0}$,
there exists $\lambda \in \R$ such that $|\lambda-\lambda_0|<\delta$,
$r_\alpha=r$, $T_\lambda|A_r \in E_r$, and
$\dist_{E_r}(T_{\lambda_0}|A_r,T_\lambda|A_r)<\delta$.

\end{thm}

Let $E_0$ be a complete metric space of continuous functions on $S^1$.  We
assume that $F_r \subset E_0$, $r>0$, and that the inclusion is continuous.

\begin{thm}

Let $r_{\lambda_0}>0$.  For every $\delta>0$,
there exists $\lambda \in \R$ such that $|\lambda-\lambda_0|<\delta$,
$r_\lambda=0$, and there exists a map $T:S^1 \to S^1$ such that $T \in E_0$,
$\dist_{E_0}(T_{\lambda_0}|S^1,T)<\delta$ and $T(e^{2 \pi i \alpha_\lambda}
z)=Q_\lambda(T(z))$.

\end{thm}

The proof of both theorems is the same as of Theorem~\ref {2.1} and we
shall not repeat it here.  We only need to replace Yoccoz's Theorem used in
(P3) by a result of Geyer \cite {G}, since Proposition 10 of Risler also
applies for Herman rings to yield (P4).  In particular, the proof works
also for the family $Q_\lambda(z)=e^{2 \pi i \lambda} z e^{a(z-1/z)}$, where
$0<|a|<1/2$ is a fixed parameter (this is the complexification of
Arnold's standard family).

{\bf Acknowledgements:} I am grateful to Xavier Buff and Jean-Christophe
Yoccoz for several suggestions.